\documentclass[]{article}
\usepackage{amsmath}
\usepackage{amsthm}
\usepackage{amstext}
\usepackage[all,cmtip]{xy}
\usepackage{amssymb}
\usepackage[justification=centering]{caption}
\newtheorem{thm}{Theorem}[subsection]

\newtheorem{lem}[thm]{Lemma}
\newtheorem{prop}[thm]{Proposition}
\theoremstyle{definition}
\newtheorem{defn}[thm]{Definition}
\theoremstyle{acknowledgement}

\theoremstyle{remark}
\newtheorem{rem}[thm]{Remark}
\theoremstyle{Conjecture}

\numberwithin{equation}{subsection}
\newtheorem{clm}[thm]{claim}
\theoremstyle{Question}

\newcommand{\vf}{\mathcal{X}^1(M)}

\begin{document}
	
	\title{ $C^1$ weak Palis conjecture for nonsingular flows}
	\small{}
	\author{Qianying Xiao\footnote{Email: xiaoqianying11@mails.ucas.ac.cn. Q. Xiao is the Communication author.  }\ \  Zuohuan Zheng\footnote{Z. Zheng acknowledges support from NNSFC($\sharp11071238$), the Key Lab. of Random Complex Structures and Data Science CAS, and National Center for Mathematics and Interdisplinary Sciences CAS.}}
	\maketitle
	\begin{center}	
		{	Academy of Mathematics and Systems Science, Chinese Academy of Sciences, Beijing, 100190, P. R. China}
	\end{center}
	
	
	\begin{abstract}
		The aim of this work is to prove $C^1$ weak Palis conjecture for nonsingular flows. Weak Palis conjecture claims that a generic vector field either is Morse-Smale or exhibits horseshoes. Central model is come up with by Crovisier to obtain homoclinic intersection for diffeomorphisms. We adapt his method for nonsingular flows.
	\end{abstract}
	\textrm{\small{\bf MCS 2010}: 37C20, 37C29, 37C27, 37D15, 37D30, 37D10}

	\section{Introduction}

	Birkhoff-Smale theorem says that transverse homoclinic intersection of periodic orbit is equivalent to the existence of horseshoe. Horseshoe is a dynamic mechanism discovered by Smale leading to deterministic chaos. On the other hand, there is a class of systems with extremely simple gradient structure: the chain recurrent set consists of finitely many critical elements, and the stable and unstable manifolds of these critical elements intersect transversely. This class is the Morse-Smale system. According to the theory of Peixoto~\cite{Pe62}, Morse-Smale systems are open and dense among $C^1$ vector fields on any closed surface. One would ask if there are typical dynamics beyond these two on general compact manifolds.
	
	\smallskip
	
	\textbf{Weak Palis Conjecture}~\cite{Pa00,Pa05}. \textit{Morse-Smale systems and systems exhibiting horseshoes are dense among all systems}.
	
	\smallskip

	For discrete dynamics, $C^1$ weak Palis conjecture has been solved by Crovisier~\cite{Cr10}. The solution of the two dimensional case is contained in Pujals-Sambarino~\cite{PS00}. Wen~\cite{Wen04} studied generic diffeomorphisms away from homoclinic tangencies and heterodimensional cycles, based on which Bonatti-Gan-Wen~\cite{BGW07} and Crovisier~\cite{Cr10} solved the three-dimensional and general dimensional cases respectively. In $C^r$ topology with $r>1$, there is no breakthrough till now.
	
	In the presence of singularities, vector fields may display robust pathological dynamics~\cite{GW79,MPP04}. For  weak Palis conjecture of flows, an advance is made by Gan-Yang~\cite{GY14} recently. Gan-Yang~\cite{GY14} prove $C^1$ weak Palis conjecture for three-dimensional singular flows. Actually, they show any nontrivial chain recurrent class of a generic flow away from homoclinic tangencies is a homoclinic class. Here homoclinic tangency is a point where the stable manifold and unstable manifold of a periodic orbit intersect nontransversely. A chain recurent class is a maximal transitive set in the weak sense, about which we will explain later.
	
	In this paper, we deal with nonsingular flows. It is believed nonsingular flows resemble diffeomorphisms~\cite{GY14}: one can use Poincar\'e maps to convert problems  to the diffeomorphic ones. And moreover we are able to adapt Crovisier's central model. Although there has been some results about higher dimensional singular flows~\cite{LGW05}, we are not able to handle singularities here.
	
	Before stating the main results, let us introduce some notations.
	
	Let $M^d$ be a $C^\infty$ compact Riemannian manifold. $\mathcal{X}^1(M)$ is the collection of $C^1$ vector fields on $M$. $\mathcal{NS}$ is the subset of nonsingular $ C^1 $ vector fields. $\mathcal{MS}$ denotes Morse-Smale systems. $\mathcal{HS}$ consists of vector fields exhibiting horseshoes.
	
	\smallskip
	\textbf{Theorem 1}.	\textit{$\mathcal{MS}\cup\mathcal{HS}$ are dense among $\mathcal{NS}$ in the $C^1$ topology}.
	\smallskip
	
	A heterodimensional cycle consists of two hyperbolic periodic orbits of different indices, say $\mathcal{O}(p)$ and $\mathcal{O}(q)$, such that $\mathcal{W}^s(\mathcal{O}(p))\cap \mathcal{W}^u(\mathcal{O}(q))\neq\emptyset$ and $\mathcal{W}^s(\mathcal{O}(q))\cap \mathcal{W}^u(\mathcal{O}(p))\neq\emptyset$. Wen~\cite{Wen02,Wen04} studied generic diffeomorphisms away from homoclinic tangencies and heterodimensional cycles.
	His results ensure that any minimally hyperbolic set admits partially hyperbolic splitting into three directions with the center one dimensional.
	
	Let $X\in \vf$, $\varphi_t$ is the flow generated by $X$ and $\Phi_t$ is the associated flow on $TM$. $ \mathrm{Sing}(X) $ is the collection of singularities of $ X $. Denote 
	\[\mathcal{N}=\mathcal{N}_X=\{v\in T_xM\mid \langle v,X(x)\rangle=0,\ x\in M\setminus \mathrm{Sing}(X)\}.\]
	 The linear Poincar\'e flow $\psi_t$ is defined as
	\[\psi_t(v)=\pi(\Phi_t(v)),\ v\in \mathcal{N},\]
	with $\pi$ the orthogonal projection to $\mathcal{N}$.

	A nonsingular set $\Lambda$ is \textit{partially hyperbolic with respect to $\psi_t$}, if there is a continuous splitting $\mathcal{N}=V^s\oplus V^c\oplus V^u$ , $T>0$ and $ 0<\lambda<1 $, such that for any $ x\in\Lambda$ and $ t\geq T$,
	\[\|\psi_t|_{V^s(x)}\|<\lambda, \ \ \  \|\psi_{-t}|_{V^u(x)}\|<\lambda,\]
	\[\|\psi_t|_{V^s(x)}\|\|\psi_{-t}|_{V^c(\varphi_t(x))}\|<\lambda,\] \[\|\psi_t|_{V^c(x)}\|\|\psi_{-t}|_{V^u(\varphi_t(x))}\|<\lambda.\]
	$V^s,\ V^c$ and $V^u$ are respectively the stable bundle, center bundle and unstable bundle.
	
	As the first step in our arguments, we obtain a result of Wen's type for nonsingular flows.
	
	\smallskip
	\textbf{Theorem 2}.
		\textit{Let $X\in \mathcal{X}^1(M)$ be a generic vector field away from homoclinic tangencies and heterodimensional cycles. Then any nonsingular nontrivial chain recurrent class either is a homoclinic class or contains a minimal set which is partially hyperbolic with one-dimensional center with respect to $\psi_t$.}
	\smallskip
	
	A chain transitive set is \textit{aperiodic} if it contains no periodic orbit. For a nonsingular aperiodic chain transitive set which is partially hyperbolic with one-dimensional center with respect to $\psi_t$, we are able to get transverse homoclinic  points in its arbitrarily small neighborhood.
	
	\smallskip
	\textbf{Theorem 3}.
		\textit{For a generic vector field $ X\in \mathcal{X}^1(M) $, any nonsingular aperiodic chain transitive set which is partially hyperbolic with one-dimensional center with respect to $\psi_t$ is contained in the closure of nontrivial homoclinic classes of $X$.} 
		
	\smallskip
	
	Our main theorem (Theorem 1) is deduced from Theorem 2 and Theorem 3.
	
	\smallskip
	\textbf{Question}.
		It is proved in ~\cite{BGY14,GY14} that any nonsingular chain transitive set of a generic three-dimensional flow away from homoclinic tangencies is hyperbolic (hence contained in a homoclinic class). One would ask whether or not any nontrivial chain recurrent class of a generic $ C^1 $ diffeomorphism away from homoclinic tangencies and heterodimensional cycles is a homoclinic class (not just approximated by homoclinic classes)? Furthermore, is any nonsingular chain recurrent class of a generic vector field away from homoclinic tangencies and heterodimensional cycles a homoclinic class?
	\smallskip
	
	This paper is organized as follows: 
	\begin{itemize}
		\item Section 2 contains the definitions, notations and a list of generic properties enjoyed by $ C^1 $ vector fields;
		\item The proof of the main theorem is given in section 3;
		\item Certain properties of a generic $C^1$ vector field away from homoclinic and heteroclinic intersections are listed in section 4. Theorem 2 is proved in this same section;
		\item In section 5, we follow the ideas of Crovisier to construct central model for nonsingular chain transitive set which is partially hyperbolic with one-dimensional center w.r.p.t. the linear Poincar\'e flow. Nontrivial homoclinic class can be obtained as in the diffeomorphic case. The proof of Theorem 3 is given at the end of this section.
	\end{itemize}

	\section{Preliminaries}
	\subsection{Definitions and notations}
	There exists $r_0>0$ such that for any $x\in M^d$,	$\exp_x:T_xM(r_0)\rightarrow M $ is injective.
	Let $X\in \mathcal{X}^1(M)$ and $x$ a regular point of $X$. For any $0<r\leq r_0$, define
	\[N_x(r)=\exp_x(\mathcal{N}_x(r)).\]
	For any $t\in \mathbb{R}$, there exists $0< \delta(x,t)\leq r_0$ such that $X$ is transverse to $N_x(\delta)$, and that the Poincar\'e map \[P_{x,t}:N_x(\delta)\rightarrow \ N_{\varphi_t(x)}(r_0)\] 
	is well-defined.
	
	\subsubsection {Poincar\'e map versus liner Poincar\'e flow}
	
	Given any $x\in M\setminus \mathrm{Sing}(X)$, there exists a small neighborhood $U$ of $x$ such that we can suppose 
	$U\subset \mathbb{R}^d$,
	$x$ the origin and that
	the Riemannian metric differs little from the Euclidean metric on $U$. Moreover, we can assume:
	\begin{itemize}
		\item  $X(0)=(1,0,\cdots,0),\ N_0(r_0)\subset \{0\}\times \mathbb{R}^{d-1}$;
		\item  For any $z\in U,\ X(z)=\sum_{i=1}^d f_i(z) e_i$ with $f_1(z)\neq0$;
		\item  There exists an affine map $A_z$ from $\mathbb{R}^{d-1}$ to $\mathbb{R}$, such that for any $ z\in U $
		\[\mathcal{N}_z=\{(A_z(v),v)\mid v\in \mathbb{R}^{d-1}\};\] 
		\item  Fixing $z=\varphi_t(0)\in U$ with $ \lvert t \rvert $ small, there exists $ r_z>0 $ such that $ N_z(r_z)\in U $.
	\end{itemize}
	
	 Since $X(z)$ is transverse to $\mathcal{N}_z$, for any $y\in \mathbb{R}^{n-1}$ with $|y|$ small, there exists $ t(y) $ close to $ t $ such that  $\varphi_{t(y)}(y)\in N_z(r_z)$, and $t(y)$ is a $C^1$ function of $y$.
	
	Define $P(y)=\varphi_{t(y)}(y)$. Then $ P(y) $ is the Poincar\'e map satisfying $P(0)=z$ and $  P(y)\in N_z(r_z)$. Futhermore,
	\[d_yP(y)=d_y\varphi_{t(y)}(y)+X(\varphi_{t(y)})\cdot d_yt(y).\]
	
	For any $v\in \mathcal{N}_0$, 
	 \[d_yP(0)(v)=d_y\varphi_t(0)(v)+X(z)\langle d_yt(0),v\rangle\in T_zN_z(r_z).\]	
	Consequently one has
	\[ d_y\varphi_t(0)(v)=d_yP(0)(v)-X(z)\langle d_yt(0),v\rangle .\]
	Note that $ T_zN_z(r_z)=\mathcal{N}_z$ and that $T_zM=\mathcal{N}_z\oplus \mathbb{R}(X(z))$. So $d_yP(0)(v)$ is the projection of $d_y\varphi_t(0)(v)$ along $X(z)$ onto $\mathcal{N}_z$.
	Since $\mathcal{N}_z$ is the perpendicular to $ X(z)$, $d_yP(0)(v)$ is the orthogonal projection of $ d_y\varphi_t(0)(v) $ to $ \mathcal{N}_z $. Therefore $ d_yP(0) $ is linear Poincar\'e map. Consequently \textit{the derivative of Poincar\'e map equals the linear Poincar\'e map}.
	
	\subsubsection{Chain recurrence~\cite{Con76}}
	
	For any $t\geq1$, $ \epsilon>0$, $ x_i\in M $ and $ t_i\geq t $ with $ 1\leq i \leq k $, $\{x_1,\cdots,x_k;t_1,\cdots, t_k\}$ is a \textit{$(t,\epsilon)$-orbit} provided 
	\[d(\varphi_{t_i}(x_i),x_{i+1})<\epsilon, \  i=1,\cdots,k-1.\]
	
	A point $y\in M $ is said to be \textit{chain attainable} from $x$ if for any $ t\geq1$ and $ \epsilon>0$, there exists an $(t,\epsilon)$-orbit $\{x_1,\cdots,x_k;t_1,\cdots, t_k\}$ such that $x_1=x,\ x_k=y$. 
	
	The \textit{chain recurrent set} of a vector field $ X $ is the collection of point that is chain attainable from itself, and is denoted by $\mathcal{R}(X)$.
	
	It is obvious that chain attainability is an equivalence relation on $\mathcal{R}(X)$. A equivalence class is called  \textit{chain recurrent class}.
	
	A compact invariant set $ \Gamma $ of $ X $ is \textit{chain transitive}, if $ y $ is chain attainable from $ x $ for any $ x $ and $ y\in \Gamma $.
	
	A compact invariant set $\Lambda$ is \textit{hyperbolic}, if there exist a continuous splitting $T_\Lambda M=E^s\oplus\langle X\rangle \oplus E^u$, $T>0 $, and $ 0<\lambda<1 $ such that for any $t>T,\ x\in \Lambda$,
	\[\|\varphi_t|_{E^s(x)}\|<\lambda,\ \|\varphi_{-t}|_{E^u(x)}\|<\lambda,\]
	with $\langle X\rangle$ the line field generated by $X$. Note that in a hyperbolic set, a singularity is isolated from regular points.
	
	A vector field $X$ is said to be hyperbolic if $\mathcal{R}(X)$ is a hyperbolic set. This is equivalent to that $X$ satisfies Axiom A and the no cycle condition: $\mathcal{R}(X)$ is decomposed into finitely many isolated basic sets which are naturally partially ordered and the system is gradient-like modulo $\mathcal{R}(X)$. Especially, a Morse-Smale system is hyperbolic.
	
	\subsubsection{Homoclinic and heteroclinic intersection}
	
	Let $\mathcal{O}(p)$ a hyperbolic orbit of $X$. The \textit{homoclinic class} of $\mathcal{O}(p)$ is the closure of transverse homoclinic orbits of $\mathcal{O}(p)$, and is denoted by $H(\mathcal{O}(p))$. By the shadowing theorem, a hyperbolic chain recurrent class is a homoclinic class.
	
	An periodic orbit $\mathcal{O}(q)$ is \textit{homoclinically related} to $\mathcal{O}(p)$, if the stable manifold of $\mathcal{O}(p)$ intersects the unstable manifold of $\mathcal{O}(q)$ transversely and the stable manifold of $\mathcal{O}(q)$ intersect transversely the unstable manifold of $\mathcal{O}(p)$. According to the $ \lambda $-Lemma, being homoclinic related determines a transitive binary relation among hyperbolic periodic orbits.
	By Birkhoff-Smale theorem, $H(\mathcal{O}(p))$ is the closure of periodic orbits homoclinically related to $\mathcal{O}(p)$.
	
	\textit{Homoclinic tangency} is a point where the stable manifold and unstable manifold of a periodic orbit intersect nontransversely. A \textit{heterodimensional cycle} consists of two periodic orbits of different indices, say $\mathcal{O}(p)$ and $\mathcal{O}(q)$, such that the stable manifold of $\mathcal{O}(p)$ intersects the unstable manifold of $\mathcal{O}(q)$ and vice versa. Due to the  dimensional inadequence, at least one heteroclinic intersection of a heterodimensional cycle is nontransversal. An arbitrarily small $C^r$ perturbation of homoclinic tangencies and heterodimensional cycles can create transverse homoclinic orbits.
	
	\subsubsection{Dominated splitting}
	
	Let $\Lambda$ be a compact invariant set, $T>0$ and $ \lambda>0 $. A continuous splitting $T_\Lambda M=E\oplus F$ is a \textit{$ (T,\lambda) $-dominated splitting} if for any $x\in\Lambda,\ t\geq T$,
	\[\|d\varphi_t|_{E(x)}\|\|d\varphi_{-t}|_{F(\varphi_t(x))}\|<e^{-\lambda t}.\]
	
	Suppose moreover $\Lambda$ is nonsingular, i.e. $\Lambda\subset M\setminus \mathrm{Sing}(X)$, and there exists a continuous splitting $\mathcal{N}_\Lambda=\Delta_{cs}\oplus\Delta_{cu}$ such that for any $ x\in\Lambda$ and $ t\geq T$,
	\[\|\psi_t|_{\Delta_{cs}(x)}\|\|\psi_{-t}|_{\Delta_{cu}(\varphi_t(x))}\|<e^{-\lambda t},\]
	then $\Lambda$ is said to admit a \textit{$ (T,\lambda) $-dominated splitting with respect to $\psi_t$}.

	\subsection{A $C^1$ generic vector field}
	
	Given a vector field $X\in\mathcal{X}^1(M)$, a \textit{critical element} of $X$ is either a singularity or a periodic orbit.
	$X$ is Kupka-Smale if all of its critical elements are hyperbolic and that the stable manifolds and unstable manifolds of the critical elements intersect transversely. According to Kupka-Smale theorem, the collection $ \mathcal{G}_{\mathcal{K-S}} $ consisting of Kupka-Smale is a dense $ G_\delta $ subset. 
	
	$\mathcal{X}^1(M)$ is a Baire space. A subset $\mathcal{R}$ of $\mathcal{X}^1(M)$ is \textit{residual} if it contains a dense $G_\delta$ subset. We say a \textit{generic} $X$ satisfies $\mathcal{P}$ if there exists a residual subset $\mathcal{G}$ such that each element of $\mathcal{G}$ satisfies $\mathcal{P}$.
	
	Since $M$ is a compact metric space, the collection of nonempty compact subsets of $M$ is nonempty and is denoted by $ \mathcal{K}(M) $. For any $K_1,\ K_2\in \mathcal{K}(M)$, define the distance $ d(K_1,\ K_2) $ as
	\[d(K_1,\ K_2)=\max\{\max_{x_1\in K_1}d(x_1,K_2),\max_{x_2\in K_2}d(K_1,x_2)\}.\]
	This metric induces the Hausdorff topology on $\mathcal{K}(M)$, and under this metric $\mathcal{K}(M)$ is compact. $\mathcal{K}(\mathcal{K}(M))$ is also a compact meric space.

	\begin{defn}
		Assume $B$ is a Baire space, $Y$ is a compact meric space. A map $\Theta:B\rightarrow \mathcal{K}(Y)$ is \textit{lower semi-continuous} at $b\in B$, if for any open set $U$ of $Y$ satisfying $U\cap\Theta(b)\neq\emptyset$, there exists a neighborhood $V$ of $B$ such that for any $z\in V$, $\Theta(z)\cap U\neq\emptyset$.
		
		$\Theta$ is said to be \textit{upper semicontinuous} at $b\in B$, if for any $K\in \mathcal{K}(Y)$ with $K\cap\Theta(b)=\emptyset$, there exists a neighborhood $V$ of $b$ such that for any $z\in V$, one has $\Theta(z)\cap K=\emptyset$.
	\end{defn}
		
		A point at which $\Theta $ is both lower and upper semi-continuous is called a continuity point. $\Theta$ is lower semi-continuous ( or upper semi-continuous ) if it is lower semi-continuous ( upper semi-continuous ) at every point of $B$. $\Theta$ is \textit{semi-continuous} if it is either lower or upper semi-continuous.

	\begin{rem}
		As indicated in~\cite{ABC11}, the collection of continuity points of a semi-continuous map is residual.
	\end{rem}
	
	We will list the $C^1$ generic properties needed in this paper. They are shared by the continuity points of certain semicontinuous dynamically related maps defined on $ \vf $.
	
	\begin{lem}
		There is a dense $G_\delta$ set $\mathcal{G}\subset\mathcal{X}^1(M)$, such that for any $ X\in \mathcal{G}$ the following statements are satisfied:
		
		\begin{enumerate}
			
			\item $X$ is Kupka-Smale;
			\item Any chain transitive set of $X$ is approximated by hyperbolic periodic orbits in the Hausdorff topology;
			\item Any chain recurrent class of $X$ containing periodic orbits is a homoclinic class;
			\item Given a compact invariant set $\Lambda\subset M$, if there is a sequence $\{X_n\}$ such that $X_n\rightarrow X$ in the $C^1$ topology with periodic orbits $\mathcal{O}(p_n)\rightarrow \Lambda$ in the Hausdorff topology, then $\Lambda$ is the Hausdorff limit of periodic orbits of $X$ itself;
			\item Suppose $\Lambda$ is a compact invariant set of $M$, there exists a sequence $\{X_n\}$ approximating $X$ in the $C^1$ topology, $H_n$ is a nontrivial homoclinic class of $X_n$ such that $\Lambda\subset\lim_{n\rightarrow \infty}H_n$. Then for any $\epsilon>0$, there exists nontrivial homoclinic class $H$ of $X$ such that $\Lambda\subset B_\epsilon(H)$.
			
		\end{enumerate}
		
	\end{lem}	
	\smallskip
	
	Item 1 comes from Kupka-Smale theorem. According to Crovisier~\cite{Cr06}, there is a dense $G_\delta$ set $\mathcal{G}_{shadow}$, such that any chain transitive set of a vector field in $\mathcal{G}_{shadow}$ can be approximated by periodic orbits in the Hausdorff topology, hence item 2 holds. Item 3 is a consequence of $C^1$ connecting lemma for pseudo-orbits~\cite{BC04}. Item 4 and 5 are by lower semicontinuity of the closure hyperbolic periodic orbits and the closure of transverse homoclinic points respectively.
	
	\section{Proof of the main theorem}
	\textit{Proof of Theorem 1}.
		Let $ \mathcal{NS} $ denotes the collection of nonsingular vector fields. $ \mathcal{HP}\subset\mathcal{X}^1(M) $ consists hyperbolic vector fields. $\mathcal{HT} $ and $ \mathcal{HC }$ are the subsets of vector fields admitting homoclinic tangencies and heteroclinicdimensional cycles respectively.
		
		Let's consider $ \mathcal{NS}\cap \mathcal{G}\setminus \overline{\mathcal{MS}\cup \mathcal{HS}} $. Since $ \overline{\mathcal{HT}}\subset\overline{\mathcal{HS}} $, $ \overline{\mathcal{HC}}\subset\overline{\mathcal{HS}} $ and $ \overline{\mathcal{HP}}\subset \overline{\mathcal{MS}\cup\mathcal{HS}} $, one has $ \mathcal{G}\setminus\overline{\mathcal{MS}\cup\mathcal{HS}}\subset\mathcal{G}\setminus\overline{\mathcal{HT}\cup\mathcal{HC}\cup\mathcal{HP}} $.

		For any $ X\in \mathcal{NS}\cap \mathcal{G}\setminus \overline{\mathcal{MS}\cup \mathcal{HS}} $, there exists a chain recurrent class $ \Lambda $ such that $ \Lambda $ is not hyperbolic. Since $ X $ is nonsingular and generic, $\Lambda$ is not reduced to a critical element. By theorem 2, either $ \Lambda $ is a homoclinic class or $\Lambda$ contains a minimal set $\Gamma$ such that $\Gamma$ is normally partially hyperbolic with one-dimensional center. The first alternative is ruled out by the assumption. Applying Theorem 3, $ \Gamma $ is contained in the closure of nontrivial homoclinic classes of $ X $, a contradiction.
		
		Therefore $ \mathcal{NS}\cap \mathcal{G}\setminus \overline{\mathcal{MS}\cup \mathcal{HS}}=\emptyset $. Since $ \mathcal{G} $ is residual, the open set $ \mathcal{NS}\setminus \overline{\mathcal{MS}\cup \mathcal{HS}}$ is empty. So Morse-Smale systems and systems admitting horseshoes are $ C^1 $ dense among nonsingular flows. The proof of Theorem 1 is finished.

	\section{Generic $C^1$ vector fields away from homoclinic and heteroclinic intersections }
	
	\subsection{Vector fields away from homoclinic tangencies}
	
	A vector field $X$ is \textit{away from homoclinic tangency} if $ X\in \mathcal{X}^1(M)\setminus \overline{\mathcal{HT}} $.
	Wen~\cite{Wen02} gives a characterization of diffeomorphisms away from homoclinic tangencies. Since Frank's lemma and the mechanism of homoclinic tangency~\cite{PS00} also hold for flows, similar assertions~\cite{GY14} holds for flows away from homoclinic tangencies.
	
	\begin{lem}~\cite[Theroem A]{Wen02}~\cite[Lemma 2.9]{GY14}
		Given any $X\in\vf\setminus\overline{\mathcal{HT}}$, there exists neighborhood $\mathcal{V}$ of $X$, $\lambda>0$, $\delta>0$, $N>1$, $T>1$, such that for any $Y\in \mathcal{V}$, any periodic point $p$ of $Y$ with period $\tau\geq T$, any partition of $[0,\tau]$
		\[0=t_0<t_1<\ldots<t_\ell=\tau,\]
		with $t_{i+1}-t_i\geq T$ for all $0\leq i\leq\ell-1$, the following statements are satisfied:
		\begin{enumerate}
			\item If $\psi_\tau(p)$ has an eigenvalue with modular in $[(1-\delta)^\tau,(1+\delta)^\tau]$, then there exists eigenspace decomposition $\mathcal{N}(p)=V^s(p)\oplus V^c(p)\oplus V^u(p)$, with $V^c(p)$ the subspace associated to the eigenvalues with modular in $[(1-\delta)^\tau,(1+\delta)^\tau]$. Moreover $\dim V^c=1 $, for any $ 0\leq i\leq \ell-1 $,
			\[\lVert\psi_{t_{i+1}-t_i}|_{\psi_{t_i}(V^s(p))}\rVert \lVert\psi_{t_i-t_{i+1}}|_{\psi_{t_{i+1}}(V^c(p))}\rVert\leq e^{-\lambda(t_{i+1}-t_i)},\]
			\[\lVert\psi_{t_{i+1}-t_i}|_{\psi_{t_i}(V^c(p))}\rVert \lVert\psi_{t_i-t_{i+1}}|_{\psi_{t_{i+1}}(V^u(p))}\rVert\leq e^{-\lambda(t_{i+1}-t_i)},\]
			and that
			\[\prod_{i=0}^{\ell-1}\|\psi_{t_{i+1}-t_i}|_{\psi_{t_i}(V^s(p))}\|\leq Ne^{-\lambda\tau},\]
			
			\[\prod_{i=0}^{\ell-1}\|\psi_{t_i-t_{i+1}}|_{\psi_{t_{i+1}}(V^u(p))}\|\leq Ne^{-\lambda\tau}.\]
			\item If $\psi_\tau(p)$ has no eigenvalue with modular in $[(1-\delta)^\tau,(1+\delta)^\tau]$, then $\mathcal{N}(p)=V^s(p)\oplus V^u(p)$ is the eigenspace decomposition, and for any $ \ i=0,\cdots,\ell-1, $
			\[\lVert\psi_{t_{i+1}-t_i}|_{\psi_{t_i}(V^s(p))}\rVert \lVert\psi_{t_i-t_{i+1}}|_{\psi_{t_{i+1}}(V^u(p))}\rVert\leq e^{-\lambda(t_{i+1}-t_i)}.\]
			Moreover one has either
			\[\prod_{i=0}^{\ell-1}\|\psi_{t_{i+1}-t_i}|_{\psi_{t_i}(V^s(p))}\|\leq Ne^{-\lambda\tau},\ \prod_{i=0}^{\ell-1}\|\psi_{t_i-t_{i+1}}|_{\psi_{t_{i+1}}(V^u(p))}\|\leq N ;\]
			or 
			\[ \prod_{i=0}^{\ell-1}\|\psi_{t_{i+1}-t_i}|_{\psi_{t_i}(V^s(p))}\|\leq N,\ \prod_{i=0}^{\ell-1}\|\psi_{t_i-t_{i+1}}|_{\psi_{t_{i+1}}(V^u(p))}\|\leq Ne^{-\lambda\tau}.\]
			
		\end{enumerate}
		
	\end{lem}
	
	\subsection{Generic nonsingular chain transitive set away from homoclinic tangencies and heterodimensional cycles}
	
	Wen~\cite{Wen04,Wen08} combined Liao's sifting lemma and Liao's shadowing lemma to an elegant result called Liao's selecting lemma. His result applies to diffeomorphisms. With no more efforts, there is a version for nonsingular flows.
	\begin{lem}~\cite[Proposition 3.7]{Wen04}
		Let $\Lambda$ be a nonsingular compact invariant set of $X$, $ \Lambda $ admits $ (T,\lambda) $ dominated splitting $\mathcal{N}_{\Lambda}= \Delta^{cs}\oplus\Delta^{cu} $ of index $i$ with respect to $\psi_t$ for some $ 1\leq i\leq d-2 $. Moreover assume:
		\begin{enumerate}
			\item there exists $ b\in\Lambda $ satisfying
			\[\prod^{n-1}_{j=0}\lVert \psi_T|_{\Delta^{cs}(\varphi_{jT}(b))}\rVert\geq1,\]
			for all $ n\geq1. $
			\item (The tilda condition) there exist $ \lambda_1 $ and $ \lambda_2 $ with $0< \lambda_2<\lambda_1<\lambda $ such that for any $ x\in \Lambda $ satisfying
			\[\prod^{n-1}_{j=0}\lVert \psi_T|_{\Delta^{cs}(\varphi_{jT}(x))}\rVert\geq e^{-\lambda_2 n T},\]
			for all $ n\geq1 $, $ \omega(x) $ contains a point $ c\in \Lambda $ satisfying
			\[\prod^{n-1}_{j=0}\lVert \psi_T|_{\Delta^{cs}(\varphi_{jT}(c))}\rVert\leq e^{-\lambda_1 n T},\]
			for all $ n\geq1. $
		\end{enumerate}
		Then for any $ \lambda_3 $ and $ \lambda_4 $ with $ 0< \lambda_4<\lambda_3<\lambda_2  $ and any neighborhood $ U $ of $ \Lambda $, there exists a hyperbolic periodic orbit $ \mathcal{O} $ of $ X $ of index $ i $ contained entirely in $ U $ with a point $ q\in \mathcal{O} $ such that
		\[\prod^{m-1}_{j=0}\lVert \psi_T|_{V^s(\varphi_{jT}(q))}\rVert\leq e^{-\lambda_4 m T},\]
		\[\prod^{[\frac{\tau(q)}{T}]-1}_{j=m-1}\lVert \psi_T|_{V^s(\varphi_{jT}(q))}\rVert\geq e^{-\lambda_3 ([\frac{\tau(q)}{T}]-m+1) T},\]
		for $ m=1,\cdots, [\frac{\tau(q)}{T}]$, with $ \tau(q) $ denotes the period of $ q $ and $ V^s $ the stable bundle of $ \mathcal{O} $.
	\end{lem}
	\begin{rem}
		With similar arguments as in ~\cite[Lemma 3.8]{Wen04}, if a nonsingular compact invariant set with dominated splitting of index $ i $ satisfies the two conditions of Liao's selecting lemma, then it intersects a homoclinic class of index $ i $.
	\end{rem}
	
	A compact invariant set $ \Lambda $ of a vector field $ X $ is \textit{minimally nonhyperbolic} if $ \Lambda $ is not hyperbolic while any proper compact invariant set of $ \Lambda $ is hyperbolic. Since being hyperbolic is an open property, by Zorn's lemma every compact invariant set which is not hyperbolic contains a minimally nonhyperbolic set.

	\subsection{Proof of Theorem 2}
		\textit{Proof of Theorem 2}. 
		Recall $\mathcal{HT}$  and $\mathcal{HC}$ denote the subset of $ C^1 $ vector fields admitting homoclinic tangencies and the subset with heterodimensional cycles respectively.
		
		Given any $X\in(\vf\setminus\overline{\mathcal{HT}\cup\mathcal{HC}})\cap\mathcal{G}$, any nonsingular chain recurrent class $\Lambda$ of $X$ that is not hyperbolic. Since $X$ is a generic  vector field, $\Lambda$ is not reduced to a periodic orbit.
		
		If $\Lambda$ contains a periodic orbit, then $\Lambda$ is a homoclinic class.
		
		Assume $\Lambda$ has no periodic points. Since $\Lambda$ is not hyperbolic, there exists a minimally nonhyperbolic set $\Gamma\subset\Lambda$: any proper compact invariant subset of $\Gamma$ is hyperbolic. $ \Gamma $ is minimal. Otherwise, by the Shadowing Lemma and the generic assumption, $\Gamma $ is contained in a homoclinic class. Hence $\Lambda$ is a homoclinic class, a contradiction.
		
		For any neighborhood $\mathcal{U}_n$ of $X$ and any neighborhood $U_n$ of $\Gamma$, there exists $X_n\in \mathcal{U}_n$ admitting a nonhyperbolic periodic orbit $\mathcal{O}(p_n)$ in $ U_n$. Otherwise, there is a neighborhood $\mathcal{U}$ of $X$ and $U$ of $\Gamma$, such that for any $Y\in\mathcal{U}$, any periodic orbit of $Y$ entirely contained in $U$ is hyperbolic. That is to say $ X $ is a star flow in $ U $~\cite{LGW05}. Consequently $\Gamma$ is hyperbolic, a contradiction.
		Let $\Gamma^\prime\subset\Gamma$ be a limit point of $\{\mathcal{O}(p_n)\}$ in the Hausdorff topology. Since $\Gamma$ is minimal, $\Gamma^\prime=\Gamma$.
		
		By the assumption that $X$ is away from homoclinic tangency, $\Gamma$ admits dominated splitting: there exists $ T>0$ and $\lambda>0 $ such that for any $x\in \Gamma$, $\mathcal{N}(x)=V_1(x)\oplus V_2(x)\oplus V_3(x)$ with dim $V_2(x)=1$, and
		\[\lVert\psi_T|_{V_1(x)}\rVert \lVert\varphi_{-T}|_{V_2(\varphi_T(x))}\rVert\leq e^{-\lambda T},\]
		\[\lVert\psi_T|_{V_2(x)}\rVert \lVert\varphi_{-T}|_{V_3(\varphi_T(x))}\rVert\leq e^{-\lambda T}.\]
		
		\begin{clm}
			$V_1$ is contracted.
		\end{clm}
		\textit{Proof of the claim.}
		Suppose $V_1$ is not contracted, then there exist $x_0\in\Gamma$ such that 
		\[\prod^{n-1}_{i=0}\lVert \psi_T|_{\psi_{iT}(V_1(x_0))}\rVert \geq1\]	
		for any $n \geq 1$.
		
		To apply Liao's Selecting Lemma, check if the tilda condition is satisfied. Suppose the tilda condition is satisfied. By Remark 4.2.2. $\Gamma$ intersects a homoclinic class, hence $\Lambda$ is a homoclinic class. So the tilda condition is not satisfied.
		For any $0<\lambda_1<\lambda$, there exists $x\in\Gamma$, such that for any $y\in\omega(x)$, there exist $n=n(y)$ with
		
		\[\prod^{n-1}_{i=0}\lVert\psi_T|_{\psi_{iT}(V_1(y))}\rVert>e^{-\lambda_1nT}.\]
		
		On the other hand, it is obvious that $ \tau(p_n)\rightarrow+\infty $, and that
		\[\prod^{[\frac{\tau(p_n)}{T}]-1}_{i=0}\lVert\psi_T^{X_n}|_{\psi_{iT}^{X_n}(V^s(p_n))}\rVert \leq CNe^{-\tau(p_n)\lambda}\]
		for $ C=\sup_{0\leq t\leq T}\lVert\psi_{-t}\rVert $.
		
		Let $ k_n\in\mathbb{N} $ such that 
		\[\prod^{k_n-1}_{i=0}\lVert\psi_T^{X_n}|_{\psi_{iT}^{X_n}(V^s(p_n))}\rVert \geq CNe^{-k_nT\lambda_1},\]
		\[\prod^{n-1}_{i=0}\lVert\psi_T^{X_n}|_{\psi_{iT}^{X_n}(V^s(p_n))}\rVert \leq CNe^{-nT\lambda_1},\ n=k_n+1,\cdots,[\frac{\tau(p_n)}{T}].\]
		
		It is easy to see that $ [\frac{\tau(p_n)}{T}]-k_n\rightarrow+\infty $ as $ n\rightarrow+\infty $. Let $ x_n=\varphi_{k_nT}(p_n) $. Suppose $\{x_n\} $ accumulate to some $x_0\in\Gamma$, then for any $ n\geq1 $
		
		\[\prod^{n-1}_{i=0}\lVert\psi_T|_{\psi_{iT}(V_1(x_0))}\rVert \leq e^{-\lambda_1nT}.\]
		Hence $ \omega(x) $ is a proper subset of $ \Gamma $, a contradiction to the fact that $ \Gamma $ is minimal.
		
		Similarly, $ V_3 $ is expanded.
		The proof of Theorem 2 is finished.
	\smallskip

	\section{Creation of homoclinc intersection}
	
	\subsection{About central model}
	Central model is come up with by Crovisier~\cite{Cr10} to study partially hyperbolic chain transitive set with one-dimensional center. It is applied to the creation of homoclinic intersections. Let's recall the related definitions and properties in~\cite{Cr10}.
	\begin{defn}~\cite[Definition 2.1]{Cr10}
		A \textit{central model} is a pair $ (\hat{K},\hat{f}) $, where $ \hat{K} $ is a compact metric space, $ \hat{f} $ is a continuous map from 
		$\hat{K}\times[0,1]  $ to $ \hat{K}\times[0,+\infty) $ such that:
		\begin{itemize}
			\item $ \hat{f}(\hat{K}\times\{0\})= \hat{K}\times\{0\}$;
			\item $ \hat{f} $ is a local homeomorphism in a small neighborhood of $ \hat{K}\times\{0\} $;
			\item $ \hat{f} $ is a skew-product: there exist two maps $ \hat{f}_1:\hat{K}\rightarrow\hat{K} $ and $ \hat{f}_2:\hat{K}\times[0,1]\rightarrow[0,+\infty) $ such that for any $ (x,t)\in \hat{K}\times[0,1]$, one has 
			\[\hat{f}(x,t)=(\hat{f}_1(x),\hat{f}_2(x,t)).\]
		\end{itemize}
	\end{defn}
	
	\begin{defn}~\cite[Definition 2.2]{Cr10}
		A central model  $ (\hat{K},\hat{f}) $ has a \textit{chain-reccurent central segment} if there exists a nontrivial segment $ I=\{x\}\times[0,a] $ contained in a chain transitive set of $ \hat{f} $.
		
	\end{defn}
	
	\begin{defn}~\cite[Definition 2.4]{Cr10}
		A subset $ \mathcal{S} $ of the  product $ \hat{K}\times[0,+\infty) $ is a \textit{strip} if for any point $ x\in \hat{K} $, the intersection $\mathcal{S} \cap(\{x\}\times[0,+\infty)) $ is an interval containing $ \{x\}\times \{0\} $.
		
	\end{defn}
	
	\begin{prop}~\cite[Proposition 2.5]{Cr10}
		Let $ (\hat{K},\hat{f}) $ be a central model with chain transitive base. Then the following two properties are equivalent:
		\begin{itemize}
			\item There is no chain-recurrent central segment;
			\item There exists strip $ \mathcal{S} $ in $\hat{K}\times[0,1]  $ that is in arbitrarily small neighborhoods of $ \hat{K}\times \{0\} $ and either $ \hat{f}(cl\mathcal{(S)})\subset int (\mathcal{S}) $ or $ \hat{f}^{-1}(cl\mathcal{(S)})\subset int (\mathcal{S}) $.
		\end{itemize}
	\end{prop}
	
	\subsection{Central model for nonsingular chain transitive set}
	
	\subsubsection{Construction of central model}
	
	Assume $ K $ is a nonsingular compact invariant set of $ X\in \vf $, $ \psi_t $ is partially hyperbolic with one dimensional center on $ K $. To be precise, assume that there exist a continuous invariant splitting $ V^s\oplus V^c\oplus V^u $ of $ \mathcal{N}_K $ and $ T>0 $ large enough, such that $ \dim V^c=1 $ and for any $ x\in K $,
	\[\lVert\psi_T|_{V^s(x)}\rVert\leq 1/4,\ \ \lVert\psi_{-T}|_{V^u(x)}\rVert\leq 1/4,\]
	moreover,
	\[\lVert\psi_T|_{V^s(x)}\rVert \lVert\psi_{-T}|_{V^c(\varphi_T(x))\oplus V^u(\varphi_T(x))}\rVert\leq 1/4,\]
	\[\lVert\psi_T|_{V^s(x)\oplus V^c(x)}\rVert \lVert\psi_{-T}|_{V^u(\varphi_T(x))}\rVert\leq 1/4.\]

	According to the theory in~\cite{HPS77}, there exists invariant plaque family associated to a dominated splitting. The construction of central model will be better illustrated by the process of cooking up  the central invariant plaque family. So let us give the details of obtaining invariant plaque families here.

		\smallskip
		
		\smallskip
		
	Since $ K\subset M\setminus \mathrm{Sing}(X) $, there exist $ r>0$ and $ \delta>0 $, such that for any $ x\in K $, the Poincar\'e map 
	\[P_{T,x}:N_x(\delta)\rightarrow N_{\varphi_T(x)}(r)\]
	is well-defined. Its lift on the normal bundle
	\[\mathcal{P}_{T,x}:\mathcal{N}_x(\delta)\rightarrow\mathcal{N}_{\varphi_T(x)}(r)\]
	is defined as $ \mathcal{P}_{T,x}= \exp_{\varphi_T(x)}^{-1}\circ P_{T,x}\circ\exp_x$. Hence we have $ \mathcal{P}_T:\mathcal{N}_K(\delta)\rightarrow \mathcal{N}_K(r)$ defined as
	\[\mathcal{P}_T(x,v)=(\varphi_T(x),\mathcal{P}_{T,x}(v)).\]
	We can take $ \delta>0 $ small enough such that $ \mathcal{P}_{T,x}-\psi_{T,x} $ is a small Lipschitz map. 
	
	Utilizing bump functions to glue $ \mathcal{P}_T $ and $ \psi_T $, we obtain $ \tilde{P}_T:\mathcal{N}_K\rightarrow \mathcal{N}_K$ such that $\tilde{P}_T-\psi_T  $ is a small Lipschitz map on each fibre and $ \tilde{P}_T|_{\mathcal{N}(\frac{1}{2}\delta)}=P_T|_{\mathcal{N}(\frac{1}{2}\delta)} $.
	
	According to~\cite[Theorem 5.5]{HPS77}, for any $ x\in K $, there exist 
	\[h_x:V^c_x\oplus V^u_x\rightarrow V^s_x,\ \  g_x:V^s_x\oplus V^c_x\rightarrow V^u_x,\]
	such that 
	\[h_x(0)=0,\ Dh_x(0)=0,\ \mathrm{lip}(h_x)<1,\]
	\[g_x(0)=0,\ Dg_x(0)=0,\ \mathrm{lip}(g_x)<1.\]
	Moreover,
	\[\tilde{P}_T(\mathrm{graph}(h_x))=\mathrm{graph}(h_{\varphi_T(x)}),\ \tilde{P}_T(\mathrm{graph}(g_x))=\mathrm{graph}(g_{\varphi_T(x)}).\]
	
	Let $ j_x:V^c_x\rightarrow V^s_x\oplus V^u_x $ such that
	\[\mathrm{graph}(j_x)=\mathrm{gragh}(h_x)\cap \mathrm{graph}(g_x).\]
	Then one has
	\[j_x(0)=0,\ Dj_x(0)=0,\ \mathrm{lip}(j_x)<1.\]
	
	According to the construction, $\{ h_x\}_{x\in K} $, $\{ g_x\}_{x\in K} $ and $\{ j_x\}_{x\in K} $ are continuous families of $ C^1 $ maps.
	
	For any $ 0<\xi<\frac{1}{2}\delta $, define
	\[W^{cu}_\xi(x)=\exp_x(\mathrm{gragh}(h_x|_{B^{cu}(\xi)})),\]
	\[W^{cs}_\xi(x)=\exp_x(\mathrm{gragh}(g_x|_{B^{cs}(\xi)})),\]
	\[W^c_x(\xi)=\exp(\mathrm{gragh}(j_x|_{B^c(\xi)})).\]
	Then one has 
	\[P_{T,x}(W^{cu}_x(\xi))\subset W^{cu}_{\varphi_T(x)}(r),\]
	\[P_{T,x}(W^{cs}_x(\xi))\subset W^{cs}_{\varphi_T(x)}(r),\]
	\[P_{T,x}(W^{c}_x(\xi))\subset W^{c}_{\varphi_T(x)}(r).\]
	
	If $ V^c_K $ is orientable, then $ \psi_T $ preserves the orientation of $ V^c_K $. Fixing an orientation of $ V^c_K $, we can define the positive half $ V^c_{K,+} $ in a natural way. There exists an isomorphism $ i:K\times[0,+\infty)\rightarrow V^c_{K,+} $ such that $ V^c_{K,+}(\frac{1}{2}\delta)=i(K\times[0,1]) $.
	
	Suppose $ V^c_K $ is not orientable. Then there is a two-fold covering map $ \hat{i}:\hat{K}\rightarrow K $ and over $ \hat{i} $ a bundle  map
	\[i:\hat{K}\times[0,+\infty)\rightarrow V^c_K.\]
	$ i $ is onto and is injective outside $ \hat{K}\times\{0\} $. We can also assume $i(\hat{K}\times[0,1])= V^c_K(\frac{1}{2}\delta) $.
	
	To simplify notations both $ K\times[0,+\infty) $ and $ \hat{K}\times[0,+\infty )$ are denoted as $ \hat{K}\times[0,+\infty )$. Let $ \mathrm{id}\times(\mathrm{id},j):V^c_K\rightarrow \mathcal{N}_K $ be defined by
	\[\mathrm{id}\times(\mathrm{id},j)(x,v)=(x,v+j_x(v)),\]
	 for any $ x\in K $ and $ v\in V^c_x $.
	Then we have the following commuting diagram:
	
	\xymatrix{
		\hat{K}\times[0,1] \ar[d]^i \ar[r]^{\hat f}              &\hat{K}\times[0,+\infty) \ar[d]^i\\
		V^c_K(\frac{1}{2}\delta) \ar[d]^{\mathrm{id}\times(\mathrm{id},j)}         &V^c_K \ar[d]^{\mathrm{id}\times(\mathrm{id},j)}\\
		\mathcal{N}_K(\frac{1}{2}\delta) \ar[d]^{\exp}           &\mathcal{N}_K \ar[d]^{\exp}\\
		M \ar[r]^{P_T}                                           &M}

	Denote $ \pi=\exp\circ (\mathrm{id}\times(\mathrm{id},j))\circ i $.
	According to the construction, $ (\hat{K},\hat{f}) $ is a central model satisfying:
	\begin{itemize}
		\item $\forall \hat{x}\in \hat{K}$, $ \pi(\{\hat{x}\}\times[0,1])\subset N_{\pi(\hat{x})}(\delta) $;
		\item $ \pi $ semiconjugates $ \hat{f}  $ and $ P_T $:
		\[\pi\circ \hat{f}|_{\{\hat{x}\}\times[0,1]}=P_{T,\pi(x)}\circ \pi|_{\{\hat{x}\}\times[0,1]};\]
		\item $ \forall \hat{x}\in \hat{K} $, $ \pi(\{\hat{x}\}\times[0,1]) $ is tangent to $ V^c $ at $ x=\pi(\hat{x},0)=\pi(\hat{x}) $. By taking $ \delta $ sufficiently small, $ \pi(\{\hat{x}\}\times[0,1]) $ is a $ C^1 $ curve almost tangent to $ V^c $.
	\end{itemize}

	\smallskip

	From now on $ X\in \mathcal{G} $ for $ \mathcal{G} $ as in section 2.2. Assume $ \Gamma $ is a compact aperiodic chain transitive set which is partially hyperbolic with one dimensional center with respect to $ \psi_t $. $ \Gamma $ is contained in a compact invariant set $ K $  which is partially hyperbolic set with one dimensional center. $ K $ can also be required to has periodic orbits arbitrarily close to $ \Gamma $ in the Hausdorff topology.
	
	As in section 5.2.1., $ K $ has central model $ (\hat{K},\hat{f}) $, the restriction of $ (\hat{K},\hat{f}) $ to $ \Gamma $ is a central model $ (\hat{\Gamma},\hat{f}) $.
	
	\subsubsection{$(\hat{\Gamma},\hat{f}) $ has chain-recurrent central segment}

	Assume $ \hat{\gamma}\subset\{\hat{x}\}\times[0,1] $ is a chain-recurrent central segment of $ (\hat{\Gamma},\hat{f}) $ for some $ \hat{x}\in \hat{\Gamma} $ and $ \gamma=\pi(\hat{\gamma}) $. Let $ \tilde{\Gamma} $ be the chain recurrent class of $ X $ containing $ \gamma\cup \Gamma $.
	
	Since $ X\in \mathcal{G} $, there exist periodic orbits arbitrarily close to $ \tilde{\Gamma} $. Take $ z\in \mathrm{int}\gamma $, for any small neighborhood $ U $ of $ z $ on $ N_{\pi(\hat{x})}(\delta) $, there exists periodic point $ p $ of $ X $ in $ U $. If $ \mathcal{N}^s_K $ is trival, then $ W^{uu}(\mathcal{O}(p)) $ intersects $ \gamma $ transversely, consequently $ \mathcal{O}(p)\subset \tilde{\Gamma} $, $ \tilde{\Gamma} $ is a homoclinic class by the generic assumption. When $ \mathcal{N}^u_K $ is trivial we have the same conclusion. 
	
	Suppose neither $ \mathcal{N}^s_K $ nor $ \mathcal{N}^u_K $ are trivial, by arguments similar to the diffeomorphic case~\cite{Cr10}, there exists $ x, y\in \gamma $ such that $ W^{uu}(\mathcal{O}(p))\cap W^{ss}(x)\neq\emptyset $, and 
	$ W^{ss} (\mathcal{O}(p))\cap W^{uu}(y)\neq\emptyset$. Therefore $ \mathcal{O}(p) $ is in the same chain recurrent class as $ \gamma $, as a result $ \Gamma $ is contained in a homoclinic class.
	
	\subsubsection{$(\hat{\Gamma},\hat{f}) $ has no chain-recurrent central segment}

	If $ (\hat{\Gamma},\hat{f}) $ has no chain-recurrent central segment, then $ (\hat{\Gamma},\hat{f}) $ or $ (\hat{\Gamma},\hat{f}^{-1}) $ has arbitrarily small trapping strip.
	
	We may as well assume $ \mathcal{S} $ is an arbitrarily small open strip and $ \hat{f}(\mathrm{cl}\mathcal{S})\subset \mathcal{S} $.
	
	\begin{lem}
		$ (\hat{K},\hat{f}) $ has arbitrarily small trapping strip.
	\end{lem}
	
	\textit{Proof.}	
	$ \mathcal{S}  $ can be looked on as the lower part of a function $ \sigma:\hat{\Gamma}\rightarrow [0,+\infty) $. $ \mathcal{S} $ being open, $ \sigma $ is lower semicontinuous. Since $ \hat{f}(\mathrm{cl}\mathcal{S}) $ and $ \hat{K}\times [0,+\infty)\setminus \mathcal{S} $ are disjoint, there are disjoint open sets $ U,V $ of $ \hat{K}\times [0,+\infty)$ such that $ \hat{f}(\mathrm{cl}\mathcal{S})\subset U $, $ \hat{K}\times [0,+\infty)\setminus \mathcal{S}\subset V $.
	Let $ \hat{K}_0 $ be the maximal invariant set within the projection of $\mathrm{cl} (V) $ on $ \hat{K} $. Define $ \rho: \hat{K}_0\rightarrow[0,+\infty)$ as
	\[\rho(x)=\inf\{t|(x,t)\in \mathrm{cl}(V), \ x\in \hat{K}_0\}.\]
		
	It is easy to see that $ \rho $ is lower semicontinous. Let $ B $ be the lower part of $ \rho $, then $ B $ is open. By taking $ V $ arbitrarliy small and that $ \hat{f} $ is continuous, $ \hat{f}(\mathrm{cl}B))\subset B $, hence $ B $ is a trapping strip of $ (\hat{K}_0,\hat{f}) $.
		
	Note that $ B $ can be arbitrarily close to $ \mathcal{S} $ by taking $ V $ arbitrarily small. Hence reducing $ \hat{K} $ if necessary,  $ \hat{K} $ admits arbitrarily small trapping strip and also denoted by $ \mathcal{S} $. The proof is finished.
		
	\smallskip
	
	\smallskip
	
	\smallskip 
	
	Define $ I=\cap_{n\geq0}\hat{f}^n(\mathcal{({S})}) $, $ \sigma_{\hat{x}}=\pi(\mathcal{S}\cap\{\hat{x}\}\times[0,1]) $, $ \gamma_{\hat{x}}=\pi(I\cap\{\hat{x}\}\times[0,1] ) $. Since $ \mathcal{S} $ is a trapping strip and $ I $ is $ \hat{f} $-invariant, one has $ P_T(\bar{\sigma}_{\hat{x}})\subset \sigma_{\hat{f}(\hat{x})} $ and $ P_T(\gamma_{\hat{x}})=\gamma_{\hat{f}(\hat{x})} $.

	It is easy to see that  $ \{\bar{\sigma}_{\hat{x}}\}_{\hat{x}\in\hat{K}} $ vary lower semicontinuous with respect to $ \hat{x} $ in the Hausdorff topology and $\{ \gamma_{\hat{x}}\}_{\hat{x}\in\hat{K}} $ vary upper semicontinous. Consequently for any $ \epsilon>0 $, there is $ \epsilon_1>0 $, such that for any $ \hat{x},\hat{y} $ with $ d(\hat{x},\hat{y})<\epsilon_1 $, $ \gamma_{\hat{x}}  $ is in the $ \epsilon $-neighborhood of $ \sigma_{\hat{y}} $.

	Suppose $ x=\pi(\hat{x})$ is a periodic point of $ X $, we have shown $ \gamma_{\hat{x}} $ is $ P_T $-invariant. Actually, $ \gamma_{\hat{x}} $ is $ P_t $ invariant for any $ t\in\mathbb{R} $.
	
	\begin{lem}
		$ P_t(\gamma_{\hat{x}})=\gamma_{\hat{\varphi_t}(x)} $ for any periodic $ x=\pi(\hat{x})\in K $, $ t\in\mathbb{R} $.
	\end{lem}

	\textit{Proof.}
		Assume $ \mathrm{Ind}(x)=\dim V^s+1 $. For any $ t\in \mathbb{R} $, $ P_t(\gamma_{\hat{x}}) $ and $ \gamma_{\hat{\varphi_t}(x) }$ are tangent to $ V^c $ at $ \varphi_t(x) $. There are neighborhoods $ U_1,\ U_2 $ of  $ \varphi_t(x) $ in $ P_t(\gamma_{\hat{x}}) $ and $ \gamma_{\varphi_t(x)}(\hat{x}) $ respectively such that 
		\[U_1\cup U_2 \subset W^s(\mathcal{O}(x))\cap N_{\varphi_t(x)}(\delta),\]
		and either 
		\[U_1\subset\bigsqcup_{y\in U_2}(W^{ss}(y)\cap N_{\varphi_t(x)}(\delta)) ,\]
		or
		\[U_2\subset\bigsqcup_{y\in U_1}(W^{ss}(y)\cap N_{\varphi_t(x)}(\delta)).\]

		Consider the negative iterate, since the length of $ \{\gamma_{\hat{z}}\}_{\hat{z}}\in \hat{K} $ are uniformly bounded, $ U_1\subset U_2 $ or $ U_2\subset U_1 $. Let $ \tilde{U}^t_1 $ be the maximal interval in $ P_t(\gamma_{\hat{x}})\cap W^s(\mathcal{O}(x)) $ and $ U^t_1 $ be the maximal interval in $ \gamma_{\hat{\varphi}_t(\hat{x})}\cap W^s(\mathcal{O}(x))$. Obviously $ \tilde{U}^t_1=P_t(U^0_1) $.
		
		\begin{clm}
			$ \tilde{U}^t_1= U^t_1$
		\end{clm}
		
		\textit{Proof of the claim}.	Otherwise, assume $ \tilde{U}^t_1\subsetneqq U^t_1 $ ($ U^t_1\subsetneqq \tilde{U}^t_1 $ is proved similary). Let $ y_0 $ be on the other end of $ \tilde{U}^t_1 $ than $ \varphi_t(x) $, then $ y_0\in W^{ss}(y) $ for some $ y\in U^t_1 $, hence $ y_0\in W^{s}(\mathcal{O}(x)) $ and there is an open interval in $ P_t(\gamma_{\hat{x}}) $ containing $ y_0 $ contained in $ W^s(\mathcal{O}(x)) $~\cite{BGW07}, a contradiction. The proof of the claim is finished.
		
		\smallskip
		
		\smallskip
		
		\smallskip 
		
		Hence $ U^t_1= \tilde{U}^t_1=P_t(U^0_1)$. Let $ p$ be the other boundry point of $ U^0_1 $ than $ x $, $ t=2\tau(x) $, since $ U^0_{2\tau}=U^{2\tau}_1=P_{2\tau}(U^0_1) $ and $ P_{2\tau} $ preserve the orientation of $ V^c(x) $, one has $ P_{2\tau}(p)=p $. Consequently $ p $ is a periodic with $ \tau(p)<3\tau(x) $ if the strip $ \mathcal{S} $ is sufficiently small. On the other hand $ \tau(x)<3\tau(p) $, hence $ \frac{1}{3} <\frac{\tau(p)}{\tau(x)}< 3$.
		By the generic assumption $ p $ is hyperbolic with $ \mathrm{Ind}(p)=\dim V^s $. 
		
		Continue this process, we may show $ \gamma_{\hat{x}} $ consists of finite hyperbolic periodic points along with the intersection of their center manifolds with $N_x(\delta)$, and $ \gamma_{\hat{x}} $ is invariant under $ P_t $ for any $ t $. The proof of the lemma is finised.
		\smallskip
		
		\smallskip
		
		\smallskip

	From the proof of Lemma 5.2.2., $ \gamma_{\hat{x}} $ contains finitely many hyperbolic periodic points. Denote the periodic point on the other end of $ \gamma_{\hat{x}} $ as $ P_{\hat{x}} $ and define stable manifolds of hyperbolic periodic points in $ \gamma_{\hat{x} }$ on $ N_{\pi(\hat{x})}(\delta)$ as $ W^s(\gamma_{\hat{x}}) $.

	\begin{lem}~\cite[Lemma 3.11]{Cr10}
		There exists $ \eta>0 $, such that for any periodic point $x=\pi(\hat{x})\in K $, the $ \eta $-neighborhood of $ \sigma_{\hat{x}} $ on $ W^{cs}_{\pi(\hat{x})}(\delta) $ is contained in $ W^s(\gamma_{\hat{x}}) $. In particular, the stable manifold of $ P_{\hat{x}} $ contains a half ball with radius $ \eta $ on $ W^{cs}_{\pi(\hat{x})}(\delta) $.
		
	\end{lem}

\textit{Proof.}
		Without loss of generality, we can assume $ V^s\perp V^c $. Since $ \mathrm{cl}(\mathcal{S})\subset\hat{f}^{-1}(\mathcal{S}) $, there exists $ \eta>0 $ such that the $ \eta $-neighborhood of $ \mathrm{cl}(\mathcal{S}) $ is contained in $ \hat{f}^{-1}(\mathcal{S}) $. Especially, for any $ \hat{x}\in \hat{K} $, the $ \eta $-neighborhood of $ \sigma_{\hat{x}}  $ in $ \pi(\{\hat{x}\}\times[0,1]) $ is contained in $ P^{-1}_{T,\pi(\hat{x})}(\sigma_{\hat{f}(\hat{x})}) $. 
		
		Assume $ z $ is in the $ \eta $-neighborhood of $ \sigma_{\hat{x}} $ on $ W^{cs}_{\pi(\hat{x})}(\delta) $, $ z_1=P_{T,\pi(\hat{x})}(z) $ lies in a small neighborhood of $ \sigma_{\hat{f}(\hat{x})} $. There exists a curve $ \gamma_1 $ starting from $ z_1 $, nearly tangent to $ V^s $, reaching $ \pi(\{\hat{f}(\hat{x})\}\times [0,1]) $ at $ z_1^\prime $. Denote $ \gamma_0=P^{-1}_{T,\pi(\hat{x})}(\gamma_1) $. It is easy to see that $ \gamma_0 $ starts at $ z $, be nearly tangent to $ V^s $, hits $ \pi(\{\hat{x}\}\times[0,1]) $ at some point $ z^\prime $. Moreover, \[z^\prime\in P^{-1}_{T,\pi(\hat{x})}(\sigma_{\hat{f}(\hat{x})}),\ \   \mathrm{length}(\gamma_1)<\frac{1}{2}\mathrm{length}(\gamma_0)\leq \frac{1}{2}\eta  .\]
		
		Let $ z_2=P^2_T(z) $. Then $ z_2 $ lies in a small neighborhood of $\sigma_{\hat{f}^2(\hat{x})}$, there is a curve $ \gamma_2 $ starting at $ z_2 $, nearly tangent to $ V^s $, intersecting $ \pi(\{\hat{f}^2(\hat{x})\}\times[0,1]) $ at $ z^\prime_2 $; $ P^{-1}_T(\gamma_2)$ goes through $ z_1 $, is nearly tangent to $ V^s $; $ P^{-2}_T(\gamma_2) $ goes through $ z $, nearly tangent to $ V^s $, intersecting $ P^{-1}_{T,\pi(\hat{x})}(\sigma_{\hat{f}(\hat{x})}) $ at some point in $ P^{-1}_{T,\pi(\hat{x})}(\sigma_{\hat{f}(\hat{x})}) $. Hence \[ z_2^\prime \in P_T(\sigma_{\hat{f}(\hat{x})}),\ \ \mathrm{length}( \gamma_2)<\frac{1}{4} \mathrm{length}(P^{-2}_T(\gamma_2))\leq \frac{1}{4}\eta .\]
		
		Continuing this process, for each $ z_n=P^n_T(z) $, there exists a curve $ \gamma_n $ nearly tangent to $ V^s $ to reach $ \pi(\{\hat{f}^n(\hat{x})\}\times[0,1]) $ at $ z_n^\prime $ satisfying 
		\[\mathrm{length}(\gamma_n)<(\frac{1}{2})^n\eta, \ \ z_n^\prime\in P^{n-1}_T(\sigma_{\hat{f}(\hat{x})}).\]
		Therefore we can see that the positive orbit of $ z $ approaches the $ \gamma_{\hat{x}} $. Consequently $ z\in W^s(\gamma_{\hat{x}}) $. The proof is done.
		\smallskip
		
		\smallskip
		
		\smallskip
	
	To conclude, in the no chain-recurrent central segment case, we have shown that for a periodic point $ x=\pi(\hat{x})\in K $,
	\begin{itemize}
		\item $ \gamma_{\hat{x}} $ is periodic and $ \gamma_{\hat{x}} $ contains finitely many periodic points;
		\item $ P_{\hat{x}} $ is a hyperbolic periodic point of index $ \dim V^s+1 $. By Lemma 5.2.4. the stable manifold of $ P_{\hat{x}} $ contains a half ball with radius $ \eta $ on 
		$ W^{cs}_{\pi(\hat{x})}(\delta) $, we denote it by $ W^s_+(P_{\hat{x}},\eta) $;
		\item $ V^u_K $ is nontrivial, otherwise the total volume of attracting basins of $ P_{\hat{x}} $'s is infinite.
	\end{itemize}
	
	\subsection{Homoclinic intersection from central model: Proof of Theorem 3}
	
	\textit{Proof of Theorem 3}.
	Let $ \Gamma $ be a nonsingular chain transitive set which is partially hyperbolic with one dimensional center with respect to the linear Poincar\'e flow.
	There exists a compact invariant set $ K $ containing $ \Gamma $ such that there exist infinitely many periodic orbits in $ K $ accumulating to $ \Gamma $ in the Hausdorff topology.
	
	According to subsection 5.2.1., $ K $ has a central model $ (\hat{K},\hat{f}) $. Its restriction to $ \Gamma $ is a central model $ (\hat{\Gamma},\hat{f}) $. If $ (\hat{\Gamma},\hat{f}) $ admits chain-recurrent central segment, then $ \Gamma $ is contained in a homoclinic class by the arguments in subsection 5.2.2.
	
	Hence we are reduced to consider the case when $ (\hat{\Gamma},\hat{f}) $ has no chain-recurrent central segment. Then $ (\hat{\Gamma},\hat{f}) $ or its inverse has arbitrary small trapping strips. We assume the former holds.
	
	Reducing $ K $ if necessary, $ K $ has sufficiently small trapping strip. 
	According to whether $ V^c_{\Gamma} $ is orientable or not, the proof will be divided into two cases. In each case, we  are able to show that $ \Gamma $ is contained in the closure of nontrivial homoclinic classes. Modulo this step, the proof of Theorem 3 is finished.
	
	\subsubsection{The nonorientable case}
	
	Suppose $ V^c_{\Gamma} $ is not orientable. For any $ x\in K $, there exist $ \hat{x}_1$ and $ \hat{x}_2\in\hat{K}$ such that $ \pi(\hat{x}_1)= \pi(\hat{x}_2)=x$. Define 
	\[\sigma_x=\sigma_{\hat{x}_1}\cup\sigma_{\hat{x}_1}, \  \gamma_x= \gamma_{\hat{x}_1}\cup\gamma_{\hat{x}_2}.\]
	
	According to the analysis at subsection 5.2.3., if $ x,y\in K $ are periodic points sufficiently close to each other, then $ \gamma_x $ is contained in a sufficiently small neighborhood of $ \sigma_y $, as a result the strong unstable manifold of any periodic point on $ \gamma_x $ intersects $ W^{cs}_y(\delta) $ at some point in the $ \eta $-neighborhood of $ \sigma_y $. By Proposition 5.2.4. the intersecting point lies in the stable manifold of some periodic point on $ \gamma_y $. With the generic assumption, the intersection is transverse.
	
	Take periodic points $ \{x_n\}^{+\infty}_{n=1}\subset K $ from different orbits such that $ x_n\rightarrow x\in \Gamma $, and that $ \tau(x_n)$tends to $+\infty $ as $ n\rightarrow +\infty $. There exist $ n_0\in\mathbb{N} $, such that for any $ m,n\geq n_0$, the unstable manifold of any periodic point on $ \gamma_{x_m} $ intersects the stable manifold of some periodic point on $ \gamma_{x_n} $ transversely and vice versa. 
	Let $ m $ be large enough such that $ \gamma_{x_{n_0}}\cap \gamma_{x_m}=\emptyset $.
	
	Define a binary relation $ \prec $ on the collection of periodic orbits of $ X $:
	\[ \mathcal{O}(p)\prec\mathcal{O}(q)\   \mathrm{if }\   W^u(\mathcal{O}(q)) \cap W^s (\mathcal{O}(p))\ \mathrm{transversely\ outside}\ \mathcal{O}(q)\cup \mathcal{O}(p). \]
	By the $ \lambda $-Lemma, the relation $ \prec $ is transitive. Let $ p\in \gamma_{x_{n_0}} $ be minimal under $ \prec $ restricted to $ \gamma_{x_{n_0}} $, i.e. if $ q\in \gamma_{x_{n_0}} $ such that $ \mathcal{O}(q)\prec\mathcal{O}(p) $, then $ \mathcal{O}(p)\prec\mathcal{O}(q) $.
	
	According to the above arguments, there exists $ q\in\gamma_{x_m} $ such that $ \mathcal{O}(q)\prec\mathcal{O}(p) $, and $ p^\prime\in\gamma_{x_{n_0}} $
	satisfying $ \mathcal{O}(p^\prime)\prec\mathcal{O}(q) $. Therefore $ \mathcal{O}(p^\prime)\prec\mathcal{O}(p) $ holds. Since $ \mathcal{O}(p) $ is minimal on $ \gamma_{x_{n_0}} $, one has $\mathcal{O}(p)\prec\mathcal{O}(p^\prime)\prec\mathcal{O}(p)$. That is to say $ \mathcal{O}(p) $ has transverse homoclinic intersection.

	\subsubsection{The orientable case}
	
	Suppose $ V^c $ is orientable, we are going to consider $ P_{\hat{x}} $'s with $ x=\pi(\hat{x}) $ a periodic point. Note that $ \tau( P_{\hat{x}})\rightarrow+\infty $ as $ \tau(x)\rightarrow+\infty $. We may as well suppose $ V^s$, $V^c$ and $ V^u $ are almost orthogonal.
	
	Since K is nonsingular, there exists $ \epsilon=\epsilon(X,K,T,\delta)>0 $ such that for any $ P_{\hat{x}} $, 
	\[F_{\hat{x}}=\{y\in M|\exists z\in N_{P_{\hat{x}}}(\epsilon),t\in[-2T,2T], y=P_{t,P_{\hat{x}}}(z)\}\]
	is a flowbox. Therefore, for any $ z\in N_{P_{\hat{x}}}(\epsilon) $, $ 0<\lvert t \rvert\leq 2T $, one has $ P_{t,P_{\hat{x}}}(z)\notin N_{P_{\hat{x}}}(\epsilon) $.
	
	For $ P_{\hat{x}}$ with $ \tau(P_{\hat{x}}) $ sufficiently large, choose an appropriate $ T_{\hat{x}}\in [T,2T) $ such that $ n_{\hat{x}}=\frac{\tau(P_{\hat{x}})}{T_{\hat{x}}} \in \mathbb{N}$. For $ i=0,\cdots,n_{\hat{x}}-1 $, define 
	\[ P_{\hat{x}}^i=\varphi_{iT_{\hat{x}}}(P_{\hat{x}}). \]
	Let $ \mathcal{O}_{\hat{x}}^i= \mathcal{O}(P_{\hat{x}})\cap N_{P_{\hat{x}}^i}(\epsilon) $. By requiring $ \tau(P_{\hat{x}}) $ large, $ \mathcal{O}_{\hat{x}}^i\neq\emptyset $ for some $ i $.
	
	For $ i=0,\cdots,n_{\hat{x}}-1 $, if $ \mathcal{O}^i_{\hat{x}}\neq\emptyset $ define $ a^i_{\hat{x}}=\min\{\mathrm{d}(P^i_{\hat{x}},z)|z\in \mathcal{O}^i_{\hat{x}}\setminus\{P^i_{\hat{x}}\}\} $, otherwise let $ a^i_{\hat{x}}=1 $.
	
	Take $ 0\leq k\leq n_{\hat{x}}-1 $ such that 
	\[ a_{\hat{x}}=a^k_{\hat{x}}=\min\{a^i_{\hat{x}}|0\leq i\leq n_{\hat{x}}-1\}=\mathrm{d}(P^k_{\hat{x}},z), \]
	for some $ z \in \mathcal{O}_{\hat{x}}^k$. It is easy to see that $ a_{\hat{x}}\rightarrow0 $ as $ \tau(P_{\hat{x}})\rightarrow+\infty $.
	
	\smallskip
	\smallskip
	
	In what follows, for any $ \Sigma\subset M $ we call the connected component of $ W^s(\mathcal{O}(x))\cap\Sigma $ containing $ x $ the stable manifold of $ x $ on $ \Sigma $.
	
	\smallskip
	
	If $ V^s $ is trivial, then $ V^u $ is codimension-one modulo the flow direction. Note that a codimension-one imbedded manifold separates the ambient manifold into two components locally. By choosing $ \hat{x} $ such that $ \tau(P_{\hat{x}}) $ is large enough, one has  
	\begin{itemize}
		\item either the unstable manifold of $ P^k_{\hat{x}} $ on $ N_{P^k_{\hat{x}}}(\epsilon) $ intersects the stable manifold of $ z $ on $ N_{P^k_{\hat{x}}}(\epsilon) $;
		\item or the unstable manifold of $ z $ on $ N_{P^k_{\hat{x}}}(\epsilon) $ intersects the stable manifold of $ P^k_{\hat{x}} $ on $ N_{P^k_{\hat{x}}}(\epsilon)$.
	\end{itemize}
	Therefore $ \mathcal{O}(P_{\hat{x}}) $ has nontrivial homoclinic orbits.
	
	\smallskip
	
	In case $ V^s $ is nontrivial, with the aid of the orientation of $ V^c $, we have the following three subcases~\cite{Cr10}:
	\begin{enumerate}
		\item the unstable manifold of $ P^k_{\hat{x}} $ on $ N_{P^k_{\hat{x}}}(\epsilon) $ intersects the stable manifold of $ z $ on $ N_{P^k_{\hat{x}}}(\epsilon) $;
		\item the unstable manifold of $ z $ on $ N_{P^k_{\hat{x}}}(\epsilon) $ intersects the stable manifold of $ P^k_{\hat{x}}$ on $ N_{P^k_{\hat{x}}}(\epsilon) $;
		\item the unstable manifold of $ P^k_{\hat{x}} $ is below the strong stable manifold of $ z $ on $ N_{P^k_{\hat{x}}}(\epsilon) $, and the unstable manifold of $ z $ is below the strong stable manifold of $ P^k_{\hat{x}} $ on $ N_{P^k_{\hat{x}}}(\epsilon) $.
	\end{enumerate}

	The first two subcases imply there exist transverse homoclinic points associated to  $ \mathcal{O}(P_{\hat{x}}) $.
	
	In the third subcase, $ z $ and $P^k_{\hat{x}}$ are in twist position~\cite{BGW07,Cr10}. By arguments similar to the diffeomorphic situation, there exists a constant $ C>0 $, a point $ A $ in the strong stable manifold of $ P^k_{\hat{x}} $ on $ N_{P^k_{\hat{x}}}(\epsilon) $, $ B $ in the unstable manifold of $ z $ on $ N_{P^k_{\hat{x}}}(\epsilon) $, and a curve $ \gamma $ starting from $ B $ to $ A $ along the positive direction of $ V^c $ such that:
	
	\begin{itemize}
		\item $ s=\mathrm{d}(P^k_{\hat{x}},A) $ and $ t=\mathrm{d}(z,B) $ satisfy
		\[\frac{1}{C}a_{\hat{x}}<s<Ca_{\hat{x}},\ \ \ \frac{1}{C}a_{\hat{x}}<t<Ca_{\hat{x}};\]
		\item Let $ \lvert \gamma \rvert $ denote the length of $ \gamma $, then $\lvert \gamma \rvert=\epsilon(a_{\hat{x}}) a_{\hat{x}}  $ with $ \epsilon(a_{\hat{x}})\rightarrow 0 $ as $ a_{\hat{x}}\rightarrow 0 $.
	\end{itemize} 
	
	\smallskip
	
	Let $ r(a_{\hat{x}})=\sqrt{\epsilon(a_{\hat{x}})} $, $ \Sigma $
	a ball on $ N_{P^k_{\hat{x}}}(\epsilon) $ centered $ A $ with radius $r_{\hat{x}}= r(a_{\hat{x}})a_{\hat{x}} $. A flowbox $ F $ is defined as
	\[F=\{y\in M|\exists y^\prime\in \Sigma,-1\leq t\leq 0, y=P_{t,P^k_{\hat{x}}}(y^\prime)\}.\]
	
	It can be verified that 
	\begin{enumerate}
		\item $ \mathcal{O}(P_{\hat{x}})\cap F=\emptyset $;
		\item $ \mathcal{O}^-(P_{-1,P^k_{\hat{x}}}(B))\cap F=\emptyset $;
		\item $ \mathcal{O}^+(A)\cap F=\emptyset $.
	\end{enumerate}
	
	Since $\lvert \gamma \rvert=\sqrt{\epsilon(a_{\hat{x}})}r_{\hat{x}}  $  and $ \sqrt{\epsilon(a_{\hat{x}})} $ tends to $ 0 $ as $ a_{\hat{x}}\rightarrow 0 $, we can apply small $ C^1 $ perturbation by $ \epsilon $-kernel lift within $ F $ to obtain a vector field $ Z $ near $ X $, such that under the dynamics of $ \varphi_t^Z $, $ A $ is on the positive orbit of $ P_{-1,P^k_{\hat{x}}}(B) $. Moreover,
	\[A\in W^s_Z(\mathcal{O}_Z(P_{\hat{x}})),\ \  P_{-1,P^k_{\hat{x}}}(B)\in W^u_Z(\mathcal{O}_Z(P_{\hat{x}})).\]
	
	Therefore $ A $ is a homoclinic point of $ \mathcal{O}_Z(P_{\hat{x}}) $. By an arbitrarily small $ C^1 $ perturbation, the homoclinic intersection is transverse.
	
	Since $\mathcal{O}_Z( P_{\hat{x}}) $ can be made arbitrarily close to $ \Gamma $ in the Hausdorff topology, $ \Gamma $ is contained in arbitrarily small neighborhood of a homoclinic class of a vector field arbitrarily close to $ X $. By the generic assumption of $ X $, $ \Gamma $ is contained in arbitrarily small neighborhood of homoclinic classes for $ X $ itself.

	\section*{Acknowledgement}
	
	We would like to express our deep gratitude to Professor Dawei Yang for posing the problem addressed in this paper and for for sharing his ideas generously. We take his advices and benefit from his critical comments. We owe Rusong Zheng for reading the main parts and Zhe Zhou for listening to the proofs in the seminar.

	
	\bibliographystyle{amsplain}
	\bibliography{referece}

\providecommand{\bysame}{\leavevmode\hbox to3em{\hrulefill}\thinspace}
\providecommand{\MR}{\relax\ifhmode\unskip\space\fi MR }
\providecommand{\MRhref}[2]{%
  \href{http://www.ams.org/mathscinet-getitem?mr=#1}{#2}
}
\providecommand{\href}[2]{#2}
\begin{thebibliography}{10}

\bibitem{ABC11}
F.~Abdenur, C.~Bonatti, and S.~Crovisier, \emph{Nonuniform hyperbolicity for
  $c^1$-generic diffeomorphisms}, Israel J. Math. \textbf{183} (2011), 1--60.

\bibitem{BC04}
C.~Bonatti and S.~Crovisier, \emph{R\'ecurrence et g\'en\'ericit\'e}, Invent.
  Math. \textbf{158} (2004), 33--104.

\bibitem{BGW07}
C.~Bonatti, S.~Gan, and L.~Wen, \emph{On the existence of non-trivial
  homoclinic classes}, Ergodic Theory Dynam. Systems \textbf{26} (2007),
  1473--1508.

\bibitem{BGY14}
C.~Bonatti, S.~Gan, and D.~Yang, \emph{Dominated chain recurrent classes with
  singularities}, Ann. Sc. Norm. Super. Pisa \textbf{accepted} (2014).

\bibitem{Con76}
C.~Conley, \emph{Isolated invariant sets and the morse index}, Regional
  Conference Series in Mathematics.(38), American Mathematical Society (1976).

\bibitem{Cr06}
S.~Crovisier, \emph{Periodic orbits and chain transitive sets of
  $c^1$-diffeomorphisms}, Publ. Math. Inst. Hautes \'Etudes Sci. \textbf{104}
  (2006), 87--141.

\bibitem{Cr10}
\bysame, \emph{Birth of homoclinic intersections: a model for the central
  dynamics of partially hyperbolic systems}, Ann. Math. \textbf{172} (2010),
  1641--1677.

\bibitem{GY14}
S.~Gan and D.~Yang, \emph{Morse-smale systems and horseshoes for three
  dimensional singular flows}, preprint (2014), ArXiv:1302.096v1.

\bibitem{GW79}
J.~Guchenheimer and R.~Williams, \emph{Structural stability of lorenz
  attractors}, Inst. Hautes \'Etudes Sci. Publ. Math. \textbf{50} (1979),
  59--72.

\bibitem{HPS77}
M.~Hirsch, C.~Pugh, and M.~Shub, \emph{Invariant manifolds}, Lecture Notes in
  Math.(583), Springer-Verlag (1977).

\bibitem{LGW05}
M.~Li, S.~Gan, and L.~Wen, \emph{Robustly transitive singular sets via approach
  of extended linear poincare flow}, Disc. Cont. Dynam. Syst. \textbf{13}
  (2005), 239--269.

\bibitem{MPP04}
C.~Morales, M.~Pacifico, and E.~Pujals, \emph{Robust transitive singular sets
  for 3-flows are partially hyperbolic attractors or repellers}, Ann. Math.
  \textbf{160} (2004), 375--432.

\bibitem{Pa00}
J.~Palis, \emph{A global view of dynamics and a conjecture on the denseness of
  finitude of attractors},  \textbf{261} (2000), 335--347.

\bibitem{Pa05}
\bysame, \emph{A global perspective for non-conservative dynamics}, Ann. I. H.
  Poincar\'e \textbf{22} (2005), 485--507.

\bibitem{Pe62}
M.~Peixoto, \emph{Structural stability on two dimensional manifolds}, Topology
  \textbf{1} (1962), 101--120.

\bibitem{PS00}
E.~Pujals and M.Sambarino, \emph{Homoclinic tangencies and hyperbolicity for
  surface diffeomorphisms}, Ann. Math. \textbf{151} (2000), 961--1023.

\bibitem{Wen02}
L.~Wen, \emph{Homoclinic tangencies and dominated splittings}, Nonlinearty
  \textbf{15} (2002), 1445--1469.

\bibitem{Wen04}
\bysame, \emph{Generic diffeomorphisms away from homoclinic tangencies and
  heterodimensional cycles}, Bull. Braz. Math. Soc.(N.S.) \textbf{35} (2004),
  419--452.

\bibitem{Wen08}
\bysame, \emph{The selecting lemma of liao}, Disc. Cont. Dynam. Syst.
  \textbf{20} (2008), 159--175.

\end{thebibliography}
\end{document}